\documentclass{hha}
\newcommand{\exte}{u}

\newcommand{\Tot}{\operatorname{Tot}}
\newcommand{\msc}{\mathcal{C}}
\newcommand{\K}{\Bbbk}
\newcommand{\Z}{\mathbb{Z}}
\newcommand{\ten}{\otimes}
\newcommand{\sk}{\operatorname{sk}}
\newcommand{\msa}{\mathcal{A}}
\newcommand{\ch}{{\operatorname{Ch}}}
\newcommand{\coch}{\operatorname{coCh}}
\newcommand{\coker}{\operatorname{coker}}
\newcommand{\ssd}{\delta}
\newcommand{\bkspace}{\mathbf{D}}
\newcommand{\e}{\mathcal{E}}
\newcommand{\Set}{\operatorname{Set}}
\newcommand{\set}[1]{\{#1\}}
\newcommand{\Top}{\operatorname{Top}}
\newcommand{\dmap}{\Theta}
\newcommand{\tp}{\otimes_\pi}
\newcommand{\setm}[2]{\{\,#1\mid#2\,\}}
\newcommand{\h}{\Lambda}
\newcommand{\ci}{\subset}
\newcommand{\im}{\operatorname{im}}
\newcommand{\tensor}{\ten}
\newcommand{\ep}{\varepsilon}
\newcommand{\voc}{\Omega}
\newcommand{\vom}{\bar{\Omega}}
\newcommand{\tvom}{\Upsilon}
\newcommand{\id}{\operatorname{id}}
\newcommand{\rfc}{\mathfrak{c}}
\newcommand{\rfb}{\mathfrak{b}}
\newcommand{\rft}{\mathfrak{t}}
\newcommand{\rfl}{\mathfrak{l}}
\newcommand{\rfm}{\mathfrak{m}}
\newcommand{\drc}{\rfc^d}
\newcommand{\drb}{\rfb^d}
\newcommand{\drt}{\rft^d}
\newcommand{\drl}{\rfl^d}
\newcommand{\rfmone}{\rfm^1}
\newcommand{\rfmtwo}{\rfm^2}
\newcommand{\wrfc}{\mathfrak{C}}
\newcommand{\wrfb}{\mathfrak{B}}
\newcommand{\wrft}{\mathfrak{T}}
\newcommand{\wdrc}{{\wrfc^d}}
\newcommand{\wdrb}{{\wrfb^d}}
\newcommand{\wdrt}{{\wrft^d}}
\newcommand{\wrfm}{\mathfrak{M}}
\newcommand{\wrfmone}{{\wrfm^1}}
\newcommand{\wrfmtwo}{{\wrfm^2}}
\newcommand{\extbrow}{\lambda}
\newcommand{\preext}{\tilde{Q}}
\newcommand{\Uinfty}{\unichar{8734}}
\newcommand{\ext}{Q}
\newcommand{\propcomparison}{5.2}
\newcommand{\prophombk}{2.3}
\newcommand{\thmsymmetric}{6.3}
\newcommand{\thmalgconv}{8.3}
\newcommand{\corvanishingofbidegrees}{9.2}
\newcommand{\propadditivity}{11.1}
\newcommand{\Sopdefinf}{12}
\newcommand{\propbringoutvert}{8.4}
\newcommand{\propeasyetwo}{4.3}
\newcommand{\thmeonebasis}{3.3}
\newcommand{\propnontrivialdiff}{9.1}
\newcommand{\lemsumsplit}{11.2}
\newcommand{\lemverticalpartial}{12.4}
\newcommand{\propcohomatvom}{6.2}
\newcommand{\thmasymmetric}{5.3}

\usepackage{tikz}
\usetikzlibrary{snakes}
\usetikzlibrary{positioning}
\usepackage[matrix,cmtip,arrow,curve]{xy}
\usepackage[unicode]{hyperref}

\begin{document}

\title{Homology operations and cosimplicial iterated loop spaces}
\author{Philip Hackney}
\email{hackney@math.ucr.edu}
\address{
		Department of Mathematics, 
		University of California, Riverside, 
		900 University Avenue,
		Riverside, CA 92521, 
		USA}
\date{\today}



\begin{abstract} If $X$ is a cosimplical $E_{n+1}$ space then $\Tot(X)$ is an $E_{n+1}$ space
and its mod 2 homology $H_*(\Tot(X))$ has Dyer-Lashof and Browder operations.
It's natural to ask if the spectral sequence converging to $H_*(\Tot(X))$
admits compatible operations.  In this paper I give a positive answer to this
question.
\end{abstract}

\received{Month Day, Year}   
\revised{Month Day, Year}    
\published{Month Day, Year}  
\submitted{Graham Ellis}  
\volumeyear{2010} 
\volumenumber{12} 
\issuenumber{2}   
\startpage{1}     
\webaddress{http://intlpress.com/HHA/v12/n2/a?}

\maketitle


\section{Introduction}

For each $n\geq 2$, let $\msc_{n+1}$ be a fixed $E_{n+1}$ operad. In this paper we extend the results of \cite{me1} to the setting of cosimplicial $\msc_{n+1}$-spaces. In particular, using the indexing convention on the homology spectral sequence so that $E_{-s,t}^2 = \pi^s H_t (X; \Z /2)$, we prove the following:

\begin{theorem}\label{T:internalthm}
Suppose that $X$ is a cosimplicial object in the category of $\msc_{n+1}$- spaces. Then there are operations in the mod-2 homology spectral sequence associated to $X$
\begin{align*}
Q^m: E_{-s,t}^r &\to E^r_{-s,m+t} & m&\in [t,t-s+n] \text{ and } n\geq s\\
Q^m: E_{-s,t}^r &\to E^w_{m-s-t,2t} & m&\in [t-s,\min(t,t-s+n)] \\
\end{align*}
where 
\[ w = \begin{cases} r & m=t-s \\ 
2r-2 & m\in [t-s+1, t-r+2] \\
r+t-m & m\in [t-r+3, t] \end{cases} \]
(noting that some of these intervals may be empty). These are homomorphisms unless $s=0$ and $m=t-s+n$. There is also a Browder operation
\[ \lambda_n: E^r_{-s,t} \ten E^r_{-s', t'} \to E^r_{-s-s', t+t'+n}\]
and, if $s=0$, then
\[ Q^{t+n} (x+y) = Q^{t+n} (x) + Q^{t+n} (y) + \lambda_n(x,y). \]
\end{theorem}

This theorem is the $E_{n+1}$ analogue of a theorem of Turner (see \cite{turner} and also our paper \cite{me1}). Prior to Turner's work, Dyer-Lashof operations were constructed in Eilenberg-Moore spectral sequences in \cite{bahri} and \cite{ligaardmadsen}, though in both cases only the `vertical' operations were given. The theory of Steenrod operations, on the other hand, has a rich literature. Appropriate starting points are \cite{dwyer} and \cite[Chapter 7]{singer}. 

As expected, the operations of Theorem~\ref{T:internalthm} are consistent for various choices of $n<  \infty$. Namely, it will be evident from our construction that, given $\msc_{n+1} \to \msc_{n+2}$, the operations are compatible with the forgetful functor
\[ (\text{$\msc_{n+2}$-spaces})^\Delta \to (\text{$\msc_{n+1}$-spaces})^\Delta. \]

In addition, we show that the operations converge to the usual Dyer-Lashof or Araki-Kudo operations in the homology of $\Tot(X)$. The precise statement is given in Theorem~\ref{T:convergencefinite}. We also show convergence of the multiplication (Theorem~\ref{T:fmult}) and of the Browder operation (Theorem~\ref{T:browderconv}). These convergence results are similar to those proved in the case $n=\infty$ in \cite{me2a}.

Much in this paper relies on its companion \cite{me1}, and the overall scheme  is similar. We work with simple \emph{Bousfield-Kan universal examples}, which we think of as cosimplicial spheres. We construct, for each $n$, \emph{external operations} for these examples and use the universal property to transport these operations into the spectral sequence for an arbitrary cosimplicial space. When that cosimplicial space is actually a cosimplicial $\msc_{n+1}$-space, we then obtain internal operations by combining the external operations with the $\msc_{n+1}(2)$-structure. 

\subsection{Notation and conventions} We generally work over the field $\K=\Z/2$, and all modules should be interpreted as being $\K$-vector spaces. 
For a chain complex $C$, we write $(\Sigma C)_{q+1} = C_q$ for the suspension, $\sk_t(C)$ for the brutal truncation with \[ \sk_t(C)_k = \begin{cases} C_k & k \leq t \\ 0 & k>t, \end{cases} \] and
 $C^v$ for the bicomplex with 
\[ (C^v)_{a,b} = \begin{cases} C_b & a=0 \\ 0 & a\neq 0. \end{cases}
\]
If $B$ is a bicomplex, we will always consider the suspension as operating vertically, so that $(\Sigma B)_{a,b+1} = B_{a,b}$.
 If $\msa$ is an abelian category, we write $N: \msa^{\Delta^{op}} \to \ch_{\geq 0}(\msa)$ for normalization and $C: \msa^\Delta \to \coch^{\geq 0} (\msa)$ for conormalization, with the convention that
\[ CY^p = \coker \Big(\bigoplus_{k=1}^p d^k: \bigoplus Y^{p-1} \to Y^p \Big). \] If $\msa$ happens to be the category $\ch$ of chain complexes over $\K$, then we interpret $C$ as landing in the category of (left-plane) bicomplexes.
If $Y$ is a cosimplicial chain complex, the indexing is given by
\[ C(Y)_{-p,q} = C(Y_q)^p. \]
Given a bicomplex $B$, we will let $TB$\label{SYM:producttot} denote the product total complex
\[ TB_m = \prod_{j} B_{j,m-j} \] together with the filtration by columns:
\[ F^{k}_m = \prod_{j\leq k} B_{j,m-j}. \] 
The spectral sequences used in this paper are derived from this filtration. In particular, the spectral sequence associated to a cosimplicial chain complex $Y$ is, by definition, the spectral sequence associated to the above filtration for the bicomplex $B=CY$.

The Bousfield-Kan universal examples (see \cite{bk} for the universal property) are cosimplicial simplicial pointed sets defined by
\[ \bkspace_{(r,s,t)} := \Sigma^{t-s} \coker \left( \sk_{s-1} \Delta_+ \to \sk_{s+r-1} \Delta_+ \right), \]
where $\Sigma$ is the Kan suspension and $\Delta_+$ is obtained by adding a disjoint basepoint to each cosimplicial degree of the standard cosimplicial simplicial set $\Delta$. 
We also define a cosimplicial chain complex
\[ D_{rst} := \Sigma^{t-s} \coker \left( \sk_{s-1} \Delta_* \hookrightarrow \sk_{s+r-1} \Delta_* \right) \cong N \K \bkspace_{(r,s,t)}, \]
where $\Delta_*$ is the normalization of the cosimplicial simplicial $\K$-module $\K\Delta$. 
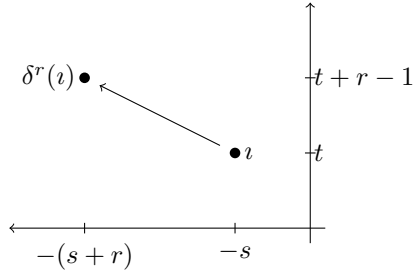
\begin{figure}[ht]
\centering
\begin{tikzpicture}
    \draw[->] (0.2,0) -- (-4,0);
    \draw[->] (0,-0.2) -- (0,3);
    \draw (-1,2pt) -- (-1,-2pt) node[below] {$-s$};
    \draw (-3,2pt) -- (-3,-2pt) node[below] {$-(s+r)$};
    \draw (2pt,1) -- (-2pt,1) node[right] {$t$};
    \draw (2pt,2) -- (-2pt,2) node[right] {$t+r-1$};
    \fill (-1,1) circle (2pt) node[right] {$\imath$};
    \fill (-3,2) circle (2pt) node[left] {$\ssd^r(\imath)$};
    \draw [->] (-1.2,1.1) -- (-2.8,1.9);
\end{tikzpicture}
\caption{Pages $2$ through $r$ of $D_{rst}$}\label{F:ue}
\end{figure}
The spectral sequence associated to $D_{rst}$ is given in Figure~\ref{F:ue}; in \cite{me1} we established the following universal property.
\begin{proposition}[Universal Property]\label{P:universalprop}
Let $Y$ be a cosimplicial chain complex and $y\in Z_{-s,t}^r(Y)$. Then there is a map \[ \dmap_y : D_{rst} \to Y \] with  $E^r (\dmap_y) (\imath) = [y]$.
\end{proposition}

Except in section~\ref{S:finiteconv}, we are not concerned about what type of spaces we use (simplicial sets or topological spaces), and consequently the mod-2 chains functor $S_*: \operatorname{Spaces} \to \ch$ will mean either the normalized simplicial chains functor $N\K$ or this functor composed with the singular functor $\Top \to \Set^{\Delta^{op}}$. The spectral sequence associated to a cosimplicial space $X$ is defined to be the spectral sequence associated to the cosimplicial chain complex $S_*(X)$.

\subsection{External operations} Since we work mod 2, the basic strategy for understanding operations follows that from \cite{may}, rather than the more complicated picture which occurs at odd primes as in \cite{cohen}. We notice that $\msc_{n+1}(2) \simeq S^n$ in the category of $\pi$-spaces (where $\pi:= \Sigma_2 = \set{e,\sigma}$). Recall the definition of $W$:
\begin{align*} W_i &= \begin{cases} e_i \cdot \K \pi & i\geq 0 \\ 0 & i <0 \end{cases}\\
d(e_i) &= (1+\sigma) e_{i-1} \end{align*}
The chain complex $W$ is a model for $S_*(\msc_\infty(2)) \simeq S_*(E\pi) \simeq S_*(S^\infty)$, and its brutal truncation $\sk_n W$ with
\[ (\sk_n W)_i = \begin{cases} W_i & i\in [0,n] \\ 0& \text{else} \end{cases} \] is a model for $S_*(S^n) \simeq S_*(\msc_{n+1}(2))$ (where all equivalences are as $\pi$-modules).
We define, for a chain complex $C$,
\[ \e^n(C) := \sk_n W \tp (C\ten C). \]
If $Z$ is a $\msc_{n+1}$-space, then, as in \cite{me1}, there is a map
\[ \e^n(S_*(Z)) \to  S_*(Z) \] induced from the structure map 
\[ \msc_{n+1}(2) \times (Z\times Z) \to Z\] and the $\pi$-homotopy equivalence $S^n \simeq \msc_{n+1}(2)$. 
The usual Araki-Kudo operations are obtained by defining chain 
maps\footnote{Notice that these are not additive until one passes to homology; in the formula one should consider all terms containing $e_k$ for $k\notin [0,n]$ to be zero.} 
(of degree $m$)
\begin{align*} q^m: C &\to \e^n(C) \\
 c &\mapsto e_{m-|c|}\otimes c \otimes c + e_{m+1-|c|} \otimes c \otimes d c 
\end{align*}
and then post-composing with the map $\e^n(C) \to C$, provided it exists.
We will call the homomorphisms $H(q^m): H_*(C) \to H_{*+m} (\e^n(C))$ \emph{external operations}.

Suppose $Y$ is a cosimplicial chain complex, and consider $\e^n(Y)$:
\[ \Delta \overset{Y}{\to} \ch \overset{\e^n}{\to} \ch. \] We regard the (homology) spectral sequence associated to $\e^n(Y)$ as the correct target for external operations. Indeed, if $X$ is a cosimplicial $\msc_{n+1}$-space and $Y= S_*(X)$, then the levelwise structure maps
\[ \msc_{n+1}(2) \times (X^p \times X^p) \to X^p \]
induce a map of cosimplicial chain complexes 
\[ \e^n(Y) \to Y. \] Thus, for this application, it is enough to construct external operations with the indicated target.

As in \cite{me1}, we construct external operations for the Bousfield-Kan universal examples $D_{rst}$. Sections \ref{S:e1e2}--\ref{S:truncfin} are devoted to computations in the spectral sequence of $\e^n(D_{rst})$. In sections \ref{S:browder} and \ref{S:truncdef} we use these calculations to define external operations for an arbitrary cosimplicial chain complex $Y$.
Finally, in section~\ref{S:finiteconv}, we extend the results of \cite{me2a} to give compatibility of our operations with those in the target.

\section{Calculation of \texorpdfstring{$E^1$ and $E^2$}{E-1 and E-2}}\label{S:e1e2}

In this section we begin the calculation of the spectral sequence for 
\[ \e^n(D_{rst}) = (\sk_n W) \tp (D_{rst} \ten D_{rst}), \] following closely the calculation for $n=\infty$ in \cite{me1}.  We assume throughout that $r\geq 2$, $n\geq 2$, and $s,t\geq 0$. 

Let us first give a basis for $E^1$. If $Y$ is a cosimplicial chain complex, then we have an isomorphism (of complexes of graded modules) $E^1(Y) = H_*C(Y)\cong CH_*(Y) $ by \cite[3.1]{bk}, where $H_*(-)$ refers to taking homology in the `chain complex direction.'  The first step towards computing $E^1(\e^n(Y))$ is thus to compute $H_*(\e^n(Y))$. By \cite[1.1]{may}, we have an isomorphism of functors out of chain complexes
\[ H_*(\e^n(-)) \cong H_*(\e^n(H_*(-))),\] so we have an isomorphism of cosimplicial graded modules
\[ H_*(\e^n(Y))  \cong  H_*(\e^n(H_*(Y))). \]
Let us collect the necessary ingredients to compute the right hand side in the case $Y=D_{rst}$.

We need the following calculation from \cite[1.3]{may}.
\begin{lemma}\label{L:mayonethree}
Let $K$ be a $\K$-module with totally ordered basis $\setm{x_j}{j\in J}$. Let $A\ci K \ten K$ have basis $\setm{x_j \ten x_j}{j\in J}$ and $B\ci K\ten K$ have basis $\setm{x_{j_1} \ten x_{j_2}}{j_1 < j_2, \text{where }j_1, j_2 \in J}$. Then
\[ H(\sk_n W\tp (K\ten K)) \cong \left( \bigoplus_{i=0}^\infty e_i \ten A \right) \oplus (e_0 \ten B) \oplus ((1+\sigma)e_n \ten B).\]
\end{lemma}
We now give a brief indication of why this lemma is true. For a $\K\pi$-module $M$ (such as $H_*(D^p) \ten H_*(D^p)$), the complex $\sk_n W\tp M$ is just 
\[ 0\to M \overset{1+\sigma}{\longrightarrow}  \cdots \to M \overset{1+\sigma}{\longrightarrow} M \overset{1+\sigma}{\longrightarrow} M \overset{1+\sigma}{\longrightarrow} M \to 0 .\] 
Thus the homology is  
\[ H_i(\sk_n W \tp M) = \begin{cases}M/(1+\sigma) &i=0 \\  \ker(1+\sigma) / \im(1+\sigma) & 0 < i < n\\ \ker(1+\sigma) & i=n \end{cases} \] and zero otherwise. 

Recall that in \cite{me1} we let  
\[ \h^p_r = \setm{\ep}{\ep:[r]\hookrightarrow [p], \ep(0)=0} \ci \Delta_p^r
\]
which was used to give a basis for $H_*(D_{rss})$. Specifically,  \cite[\prophombk]{me1} states that (for $s\geq 0$ and $r\geq 2$) the homology of $D_{rss}^p$ is given by
\[ H_k(D_{rss}^p) \cong \begin{cases}
\K  \h_{s+r}^p & k=s+r-1 \\
\K \h_s^p & k=s \\
0 & \text{else.}
 \end{cases} \]
Since Lemma~\ref{L:mayonethree} mentions an \emph{ordered} basis, we might, for concreteness, define an order at this point.
We associate to $\ep: [m]\hookrightarrow [p]$ the word of length $p+1$, whose $i^{\text{th}}$ letter is $0$ if $i\notin \im \ep$ and $1$ if $i\in \im \ep$. 
For a fixed $m$ we declare the order on injections to be given by the reverse lexicographic order on their associated words. We then define an order on the basis $\h_s^p \sqcup \h_{s+r}^p$ for $H_*(D^p)$ by insisting that $\h_s^p < \h_{s+r}^p$.

\begin{theorem}\label{Ts:e1basis}
The $E^1$ page of the spectral sequence associated to $\e^n(D_{rst})$ has a basis consisting of bigraded sets (column, bottom, top, middle) \[ \wrfc, \wdrc, \wrfb, \wdrb, \wrft, \wdrt, \wrfmone, \text{ and } \wrfmtwo.\] 
 We give an exhaustive list of their elements.
For $m\in [0,n]$, we have
\[ e_m \tensor \id_{[s]} \tensor \id_{[s]} \in \wrfc_{-s,2t+m} \]
\[ e_m \tensor  \id_{[s+r]} \tensor  \id_{[s+r]} \in (\wdrc)_{-s-r,2t+2r+m-2}. \]

For each pair $\ep < \ep' \in \h_s^p$ and $[p]=\im \ep \cup \im \ep'$ we have
\[ e_0 \tensor \ep \tensor \ep' \in \wrfb_{-p,2t} \]
\[ (1+\sigma) e_n \ten \ep \ten \ep' \in \wrft_{-p, 2t+n}.\]
Here $-p$ is between $-s-1$ and $-2s$.

For $\ep\in \h_s^p, \gamma\in \h_{s+r}^p$ and $[p]=\im \ep \cup \im \gamma$, we have
\[ e_0 \tensor \ep \tensor  \gamma \in (\wrfmone)_{-p,2t+r-1} \]
\[ (1+\sigma)e_n \tensor \ep \tensor  \gamma \in (\wrfmtwo)_{-p,2t+r+n-1} \] which live in degrees with $-p$ between $-s-r$ and $-2s-r$.

For $\gamma <\gamma' \in \h_{s+r}^p$ and $[p]=\im \gamma \cup \im \gamma'$ we have
\[e_0 \tensor  \gamma \tensor  \gamma' \in (\wdrb)_{-p,2t+2r-2} \]
\[(1+\sigma)e_n \tensor  \gamma \tensor  \gamma' \in (\wdrt)_{-p,2t+2r+n-2} \]
with $-p$ between $-s-r-1$ and $-2s-2r$. If $r=\infty$ then the basis only consists of $\wrfc$, $\wrfb$, and $\wrft$.
\end{theorem}
\begin{proof} Applying Lemma~\ref{L:mayonethree}, the homology of $\e^n(H_*(D_{rst}^p))$ has a basis given by the disjoint union of the following sets:
\begin{gather*} 
\setm{e_m\tensor \ep \tensor \ep}{\ep\in \h_s^p, m\in [0,n]} \\
\setm{e_m \tensor \gamma \tensor  \gamma}{\gamma\in \h_{s+r}^p, m\in [0,n]} \\
\setm{e_0 \tensor \ep \tensor \ep'}{\ep, \ep'\in \h_s^p, \ep < \ep'}\\
\setm{e_0 \tensor \ep \tensor \gamma}{\ep\in \h_s^p, \gamma\in \h_{s+r}^p}\\
\setm{e_0 \tensor \gamma \tensor \gamma'}{\gamma, \gamma'\in\h_{s+r}^p, \gamma < \gamma'}\\
\setm{(1+\sigma)e_n \tensor \ep \tensor \ep'}{\ep, \ep'\in \h_s^p, \ep < \ep'}\\
\setm{(1+\sigma)e_n \tensor \ep \tensor \gamma}{\ep\in \h_s^p, \gamma\in \h_{s+r}^p}\\
\setm{(1+\sigma)e_n \tensor \gamma \tensor \gamma'}{\gamma, \gamma'\in\h_{s+r}^p, \gamma < \gamma'}.
\end{gather*}
Now conormalize $H_*(\e^n(D_{rst}))$ as in \cite{me1}.
\end{proof}

It is helpful to visualize this page. 
There are five potential pictures for $E^1$ that can be drawn, depending on the relationship between $n$ and $r$.  One such picture  (for the case $n>2r-2$)  is given in Figure~\ref{F:e1trunc}, where we have indicated modules with rank greater than zero by snaky lines and modules of rank equal to one by straight lines.

\begin{figure}[ht]
\centering
{
\begin{tikzpicture}
    \draw[->] (-.65,0) -- (-8,0);   
    \draw (0.2,0) -- (-.35,0);      
    \draw (-.75,-.1) -- (-.55,.1);  
    \draw (-.45,-.1) -- (-.25,.1);  
    \draw[->] (0,-0.2) -- (0,5);
    \draw (-1,2pt) -- (-1,-2pt) node[below=.01] {$-s$};
    \draw (-2.5,2pt) -- (-2.5,-2pt) node[below=.01] {$-s-r$};
    \draw (-4,2pt) -- (-4,-2pt) node[below=.01] {$-2s$};
    \draw (-5.5,2pt) -- (-5.5,-2pt) node[below=.01] {$-2s-r$};
    \draw (-7,2pt) -- (-7,-2pt) node[below left=.01 and -.6] {$-2s-2r$};
    \draw (2pt,1) -- (-2pt,1) node[right] {$2t$};
    \draw (2pt,2) -- (-2pt,2) node[right] {$2t+r-1$};
    \draw (2pt,3) -- (-2pt,3) node[right] {$2t+2r-2$};
    \draw (2pt,3.6) -- (-2pt,3.6) node[right] {$2t+n$};
    \draw [snake=snake,very thick] (-5.5,2) -- (-2.5,2);
    \draw [snake=snake,very thick] (-5.5,4.6) -- (-2.5,4.6);
    \draw [snake=snake,very thick] (-4,1) -- (-1,1) [snake=none] -- (-1,3.6);
    \draw [snake=snake,very thick] (-7,3) -- (-2.5,3) [snake=none]--  (-2.5,5.6);
    \draw [snake=snake,very thick] (-7,5.6) -- (-2.5,5.6);
    \draw [snake=snake,very thick] (-4,3.6) -- (-1,3.6);
\end{tikzpicture} 
}
\caption{$E^1(  \e^n(D_{rst}))$}\label{F:e1trunc}
\end{figure}

The basis elements with $e_0$ and $e_n$ give us up to six interesting rows to talk about. There may be overlaps of snaky segments when $n=r-1$ or $2r-2$, but this does not affect the calculation of the homology of $\ssd^1$. 

It turns out that we have already computed the homology of each of the six snaky line segments in \cite{me1}, and we now recall three families of complexes from that paper. First, we have
\[ \voc_{s,s'} = C(H_s(D_{\infty ss}) \ten H_{s'}(D_{\infty s's'})). \] When $s=s'$ there is an action of $\pi$, and we then have the associated complexes
\[ \vom_s = \voc_{s,s} / \pi \]
and
\[ \tvom_s = \ker(1+\sigma: \voc_{s,s} \to \voc_{s,s}). \]

The following proposition is a conglomeration of \cite[\thmasymmetric, \propcohomatvom, \thmsymmetric]{me1}.
\begin{proposition}\label{P:cohomologies}
Fix $s,s'\geq 0$. The cohomologies of $\voc_{s,s'}$, $\vom_s$, and $\tvom_s$ are given by 
\begin{align*} H^p \voc_{s,s'} &= \begin{cases} \K & \text{if } p = s+s', \\ 0 & \text{otherwise}\end{cases} &
 H^p \vom_s &= \begin{cases} \K & p \in [s,2s] \\ 0 &\text{otherwise} \end{cases} \\
H^p (\tvom) &= \begin{cases} \K & \text{for } p\in [s+2,2s] , s > 1, \\ \K & \text{for } p=0=s, \\ 0 & \text{else.}\end{cases}    
\end{align*}
\end{proposition}

Using our current notation, \cite[\propcomparison]{me1} asserts that we have isomorphisms of complexes
\begin{equation}
\begin{aligned}
\vom_s &\to \left( \K(\wrfb \sqcup \wrfc_{-s,2t}), \ssd^1 \right) \\
\voc_{s,s+r} &\to \left( \K(\wrfmone), \ssd^1 \right) \\
\vom_{s+r} &\to \left( \K(\wdrb \sqcup (\wdrc)_{-s-r,2t+2r-2}), \ssd^1 \right). \\
\end{aligned} \label{E:propcomparison}
\end{equation}
This takes care of the identification of three of the six snaky line segments. For the other three, we have
\begin{proposition}\label{Ps:comparison}
There are isomorphisms of complexes
\begin{align*}
\tvom_s &\to \left( \K(\wrft \sqcup \wrfc_{-s,2t+n}), \ssd^1 \right) \\
\voc_{s,s+r} &\to \left( \K(\wrfmtwo), \ssd^1 \right) \\
\tvom_{s+r} &\to \left( \K(\wdrt \sqcup (\wdrc)_{-s-r,2t+2r+n-2}), \ssd^1 \right). \\
\end{align*}
\end{proposition}
\begin{proof}
The map 
\begin{align*}
H_s(D_{rss}) \ten H_{s+r-1} (D_{rss}) &\to H_* (\e^n(H_*(D_{rss})))\\
\ep \ten \gamma' & \mapsto (1+\sigma) e_n \ten \ep \ten \gamma'
\end{align*}
induces the middle isomorphism, as in the proof of \cite[\propcomparison]{me1}.
Furthermore, if $M^\bullet$ is one of the cosimplicial modules 
\[ H_s(D_{rss}) \ten H_s(D_{rss})  \text{ or } H_{s+r-1}(D_{rss}) \ten H_{s+r-1}(D_{rss}), \] then we have a cosimplicial map
\begin{align*}
\ker (1+\sigma: M\to M) &\to H_* (\e^n(H_*(D_{rss}))) \\
\zeta \ten \zeta' & \mapsto e_n \ten \zeta \ten \zeta'.
\end{align*}
This is an inclusion by Lemma~\ref{L:mayonethree}, and remains so after conormalizing. Finally, it is easy to see that the conormalized map
\[ \tvom \to E^1(\e^n(H_*(D_{rss}))) \] has the appropriate image.
\end{proof}

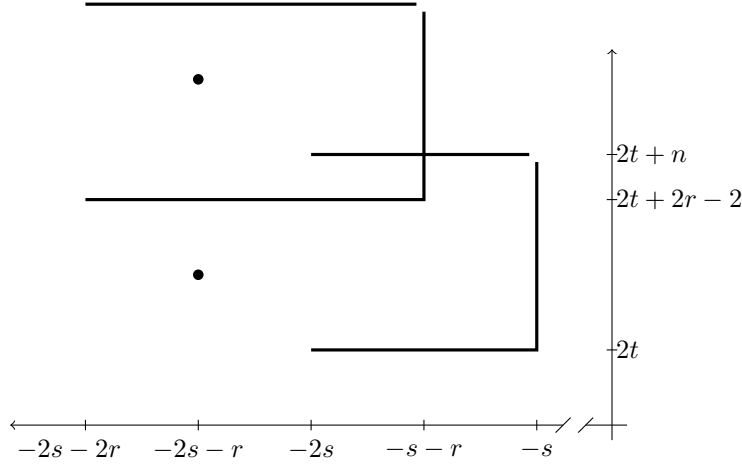
\begin{figure}[ht]
\centering
{
\begin{tikzpicture}
    \draw[->] (-.65,0) -- (-8,0);   
    \draw (0.2,0) -- (-.35,0);      
    \draw (-.75,-.1) -- (-.55,.1);  
    \draw (-.45,-.1) -- (-.25,.1);  
    \draw[->] (0,-0.2) -- (0,5);
    \draw (-1,2pt) -- (-1,-2pt) node[below=.01] {$-s$};
    \draw (-2.5,2pt) -- (-2.5,-2pt) node[below=.01] {$-s-r$};
    \draw (-4,2pt) -- (-4,-2pt) node[below=.01] {$-2s$};
    \draw (-5.5,2pt) -- (-5.5,-2pt) node[below=.01] {$-2s-r$};
    \draw (-7,2pt) -- (-7,-2pt) node[below left=.01 and -.6] {$-2s-2r$};
    \draw (2pt,1) -- (-2pt,1) node[right] {$2t$};
    \draw (2pt,3) -- (-2pt,3) node[right] {$2t+2r-2$};
    \draw (2pt,3.6) -- (-2pt,3.6) node[right] {$2t+n$};
    \fill (-5.5,2) circle (2pt);
    \fill (-5.5,4.6) circle (2pt);
    \draw [very thick] (-4,1) -- (-1,1) -- (-1,3.5);
    \draw [very thick] (-7,3) -- (-2.5,3) -- (-2.5,5.5);
    \draw [very thick] (-7,5.6) -- (-2.6,5.6);
    \draw [very thick] (-4,3.6) -- (-1.1,3.6);
\end{tikzpicture} 
}
\caption{$E^2( \e^n(D_{rst}))$ for $n>2r-2$}\label{F:e2trunc}
\end{figure}

All of the work to compute $E^2$ has now been completed. Figure~\ref{F:e1trunc} becomes Figure~\ref{F:e2trunc}.


\begin{theorem}\label{Ts:e2page} Suppose $s>0$. Then a basis for $E^2 (\e^n( D_{rst}))$ is given by the disjoint union of bigraded sets $\rfc,\drc, \rfb,\drb, \rft,\drt, \rfmone,$ and $\rfmtwo$,  which consist of a  single element in each of the indicated bidegrees:
\begin{align*}
\rfc:& \set{-s} \times [2t,2t+n-1] & \drc:& \set{-s-r} \times [2t+2r-2, 2t+2r+n-3] \\
\rfb:& [-2s,-s-1] \times \set{2t} & \drb:& [-2s-2r,-s-r-1] \times \set{2t+2r-2} \\
 \rft:& [-2s,-s-2] \times \set{2t+n}  & \drt:& [-2s-2r,-s-r-2] \times \set{2t+2r+n-2} \\
 \rfmone:& \set{(-2s-r,2t+r-1)} 
  & \rfmtwo:& \set{(-2s-r, 2t+n+r-1)} \\
\end{align*}
If $s=0$ then the statement is the same except that $\rfc$ is also nonempty in bidegree $(-s,2t+n)$ (and of course  $\rfb=\emptyset = \rft$). If $r=\infty$ then the basis is given by $\rfc$, $\rfb$, and $\rft$.
\end{theorem}
\begin{proof}
Ignoring the vertical grading, $(E^1, \ssd^1)$ is the direct sum of the complexes \begin{align*} 
&\left( \K(\wrfb \sqcup \wrfc_{-s,2t}), \ssd^1 \right) &&  \left( \K(\wdrb \sqcup (\wdrc)_{-s-r,2t+2r-2}), \ssd^1 \right) \\ 
&\left( \K(\wrft \sqcup \wrfc_{-s,2t+n}), \ssd^1 \right) && \left( \K(\wdrt \sqcup (\wdrc)_{-s-r,2t+2r+n-2}), \ssd^1 \right) \\
&\bigoplus_{i=1}^{n-1} (\K\wrfc_{-s,2t+i}, 0) && \bigoplus_{i=1}^{n-1} (\K\wdrc_{-s-r,2t+2r-2+i}, 0) \\
&\left( \K\wrfmone, \ssd^1 \right) && \left( \K\wrfmtwo, \ssd^1 \right).
\end{align*}
The result now follows from 
Propositions~\ref{P:cohomologies}, \ref{Ps:comparison}, and \eqref{E:propcomparison}.
\end{proof}


A variation of the proof of \cite[8.6]{me1} gives that
\begin{proposition}\label{Ps:einfzero} If  $r$ is finite then
$TE^\infty(\e^n(D_{rst})) = HTC(\e^n(D_{rst})) = 0$. \end{proposition}
In \cite{me1} (the case $n=\infty$), the structure of $E^2$, combined with the vanishing of $E^\infty$, determined the higher differentials in the spectral sequence.  This is not the case here. Considering naturality would allow us to compute the higher differentials, but even a statement along the lines of \cite[\propnontrivialdiff]{me1} is horrendously complicated with many cases. Luckily, we will only need partial information about the differentials in this spectral sequence. We make some computations in section~\ref{S:truncfin}, but before that we restrict to the case $r=\infty$ where we \emph{are} able to give a complete answer. Section~\ref{S:truncinf} is also an important component of the convergence results in section~\ref{S:finiteconv}.



\section{Special case: \texorpdfstring{$r=\infty$}{r=\Uinfty}}\label{S:truncinf}
\begin{figure}[ht]
\centering
{
\begin{tikzpicture}
    \draw[->] (0.2,0) -- (-8,0);
    \draw[->] (0,-0.2) -- (0,5);
    \draw (-1,2pt) -- (-1,-2pt) node[below] {$-s$};
    \draw (-4,2pt) -- (-4,-2pt) node[below] {$-2s$};
    \draw (2pt,1) -- (-2pt,1) node[right] {$2t$};
    \draw (2pt,3.6) -- (-2pt,3.6) node[right] {$2t+n$};
    \draw [very thick] (-4,1) -- (-1,1) -- (-1,3.5);
    \draw [very thick] (-4,3.6) -- (-1.1,3.6);
\end{tikzpicture} 
}
\caption{$E^2( \e^n(D_{\infty st}))$}\label{F:e2truncinf}
\end{figure}
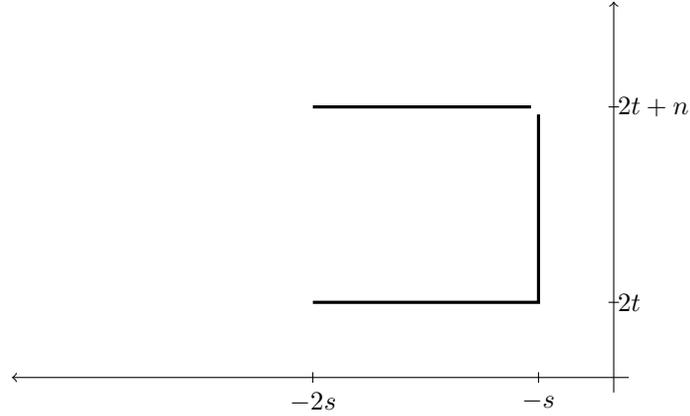

In Figure~\ref{F:e2truncinf} we visualize the statement of Theorem~\ref{Ts:e2page} when $r=\infty$. Notice that $E^2\neq E^\infty$ when $n<\infty$ and $s> 0$.
The reason for this is that there are $n+1$ nonzero terms on $E^\infty$, but $E^2$ has $2s+n-2$ nonzero classes. To justify the first part of this sentence, notice that an easy variant of the proof of \cite[\thmalgconv]{me1} gives 
\begin{theorem}\label{Ts:algconv}
Let $2\leq r \leq \infty$ and $0\leq n \leq \infty$. Then
\[ H_* TC(\e^n(D_{rst})) \cong TE^\infty (\e^n(D_{rst})). \]
\end{theorem}

Let $\rfl = \rfb \sqcup \rfc$ denote the generating set for the lower right portion of $E^2(\e^n(D_{\infty s t}))$. Since there is some overlap in total degrees, we need to check that the terms that survive to $E^\infty$ are all in $\rfl$.

\begin{theorem}\label{T:iinfty} Consider the spectral sequence for $\e^n(D_{\infty s t})$. The classes from $E^2$ which survive and are nonzero at $E^\infty$ are exactly those of $\rfl$ in total degrees \[ [2t-2s, 2t-2s+n].\]
\end{theorem}
\begin{proof}
By \cite[\propbringoutvert]{me1} we know that $(\sk_n W)^v \tp C(D\ten D) \cong C(\sk_n W \tp (D\ten D))$, so according to Theorem~\ref{Ts:algconv} it is enough to compute the cohomology of 
\[ T((\sk_n W)^v \tp C(D\ten D)), \]  where $D=D_{\infty ss}$.
If we set
\[ b_m^n = \dim H_m T((\sk_n W)^v \tp C(D\ten D)), \] we already know from Theorem~\ref{Ts:algconv}  
that
\[ b_m^\infty = \begin{cases} 1 & m\geq 0 \\ 0 & m < 0.  \end{cases} \] 
For $r\geq 2$ we have, by \cite[\propeasyetwo]{me1}, that \[ E^r(D) \otimes E^r(D) \cong E^r(D\otimes D) = E^r(\sk_0\tp (D\otimes D)).\] Thus, Theorem~\ref{Ts:algconv} gives that
\[ b_m^0 = \begin{cases} 1 & m=0 \\ 0 & m\neq 0. \end{cases} \]
We will now use comparison and induction to interpolate between these two extremes and show that
\[ b_m^n = \begin{cases} 1 & m\in [0, n] \\ 0 & \text{otherwise.}\end{cases} \] 

We have two exact sequences of complexes
\begin{equation}\label{E:seqinf} 0 \to \sk_n W \to W \to \Sigma^{n+1} W \to 0\end{equation}
and 
\begin{equation}\label{E:seqfin} 0 \to \sk_0 W \to \sk_n W \to \Sigma \sk_{n-1} W \to 0 \end{equation} where the first map is the obvious inclusion in both cases.
Freeness of $W_i$ over $\K \pi$ tells us that if we apply $(-)^v \tp C(D \ten D)$ to either of these exact sequences we will still have a short exact sequence of bicomplexes. Furthermore, taking products is exact, so applying $T$ we again get short exact sequences of complexes -- write \[ B^n = T((\sk_n W)^v \tp C(D_{\infty ss} \ten D_{\infty ss})). \] 

We first apply $T((-)^v \tp C(D \ten D))$ to SES(\ref{E:seqfin}). 
The short exact sequence of complexes
\[ 
0 \to B^0 \to B^n \to \Sigma B^{n-1} \to 0
\]
gives us a long exact sequence in homology, so we have
\[ \xymatrix@R=.8pc{
H_m B^0 \ar@{->}[r] \ar@{=}[d] & H_m B^n \ar@{->}[r] & H_m \Sigma B^{n-1} \ar@{->}[r] \ar@{=}[d] & H_{m-1} B^0 \ar@{=}[d] \\
0 &  & H_{m-1} B^{n-1} & 0
}\] for $m-1 > 0$. If we assume inductively that 
\[ b_m^{n-1} = \begin{cases} 1 & m\in [0,n-1] \\ 0 & \text{otherwise} \end{cases} \] then we see that $b_m^n$ is zero for $m\in [n+1, \infty)$ and one for $m\in [2,n]$.

Next we apply $T((-)^v \tp C(D \ten D))$ to SES(\ref{E:seqinf}), so we have
\[ 0 \to B^n \to B^\infty \to \Sigma^{n+1} B^\infty \to 0 \]
and in the long exact sequence in homology
we get
\[ \xymatrix@R=.8pc{
H_{m+1} \Sigma^{n+1} B^\infty \ar@{->}[r] \ar@{=}[d] & H_m B^n \ar@{->}[r]& H_m B^\infty \ar@{->}[r]& H_m \Sigma^{n+1} B^\infty \ar@{=}[d] \\
H_{m-n} B^\infty &  & & H_{m-n-1} B^\infty
}\]
This tells us that $b_m^n = b_m^\infty$ for $m < n$, since there $b_{m-n}^\infty = 0 = b_{m-n-1}^\infty$. 

We thus know $b_m^n$ for $m\in [2,\infty) \cup [0, n-1]$, so for $n\geq 2$ we know it for all $m$. The only thing we are missing is $b_1^1$, but this follows from the exact sequence 
\[ \xymatrix@R=.8pc{ 
0 \ar@{->}[r] & H_1 B^1 \ar@{->}[r] & H_1 \Sigma B^0 \ar@{->}[r] \ar@{=}[d] & H_0 B^0 \ar@{->}[r] \ar@{=}[d] & H_0 B^1 \ar@{->}[r] \ar@{=}[d] & 0 \\
& & \K & \K & \K 
}\]
\end{proof}


The following corollary says that there is only one way for something in $\rft$ to be hit by something in $\rfl$. 
Recall that $\rfl = \rfb \sqcup \rfc$ and $\rft$ consist of at most one element in each bidegree. 

\begin{corollary}\label{C:diffinfty} Consider the spectral sequence for $\e^n(D_{\infty st})$. 
If \[ a \in [\max (n+1-s, 0), n-1],\] then examination of bidegrees tells us that
\[ \rfc_{-s,2t+a} = \set{c} \qquad \text{and} \qquad \rft_{a-s-n-1,2t+n}=\set{v} \] are nonempty.
The element $c$ survives to $E^{n+1-a}$ and
\[ \ssd^{n+1-a}[c] = [v] \neq 0. \] 

If $a\in [n+1,s-1]$, then
\[ \rfb_{a-2s, 2t} = \set{c} \qquad \text{and} \qquad \rft_{a-2s-n-1, 2t+n} = \set{v} \] are nonempty. The element $c$ survives to $E^{n+1}$ and
\[ \ssd^{n+1}[c] = [v] \neq 0. \]

\end{corollary}
\begin{proof} When $s=0$ this corollary says that there are no nontrivial differentials, which is obvious since $E^2$ consists of a single column.

Assume $s>0$ and $t=0$. First note that $\rfl$ lives in total degrees $[-2s,-s+n-1]$ and $\rft$ lives in total degrees $[-2s+n, -s+n-2]$. All differentials out of $\rft$ are zero, so those elements in total degrees $(-2s+n,-s+n-2]$ must be hit by something (Theorem~\ref{T:iinfty}). Thus elements of $\rfl$ in total degrees 
$ (-2s+n+1, -s+n-1]$ have some nonzero differential. Collectively, these classes hit the classes of elements of $\rft$ in total degrees $(-2s+n, -s+n-2]$.
The classes that are unaccounted for are the elements of $\rfl$ in total degrees $[-2s,-2s+n+1]$ and the element of $\rft$ in total degree $-2s+n$. The element $c$ which lives in total degree $-2s+n+1$ of $\rfl$ doesn't survive to $E^\infty$, and since it lives in the second page it cannot map to something in $\rfl$. Hence there is a nontrivial differential
\[ 0 \neq \ssd^*([c]) = [v] \]
where $v\in \rft$ is in total degree $-2s+n$. The statement then follows by passing from total degree to bidegree.
\end{proof}

\section{Remaining pages, case \texorpdfstring{$r<\infty$}{r<\Uinfty}}\label{S:truncfin}

We now turn our attention to partial external operations. For a cosimplicial chain complex $Y$, these are operations whose target is the spectral sequence for
\[ \e^n(Y) = \sk_n W \tp (Y \ten Y).\] They are of particular interest when we have a map
\[ \e^n (Y) \to Y,\]
such as when $Y=S_*X$, where $X$ is a cosimplicial $E_{n+1}$-space.

\begin{figure}
\centering
\scalebox{.7}
{
\begin{tikzpicture}
	\draw[->] (-.65,0) -- (-8,0);		
	\draw (0.2,0) -- (-.35,0);		
	\draw (-.75,-.1) -- (-.55,.1);	
	\draw (-.45,-.1) -- (-.25,.1);	
    \draw[->] (0,-0.2) -- (0,5);
    \draw (-1,2pt) -- (-1,-2pt) node[below=.01] {$-s$};
    \draw (-2.5,2pt) -- (-2.5,-2pt) node[below=.01] {$-s-r$};
    \draw (-4,2pt) -- (-4,-2pt) node[below=.01] {$-2s$};
    \draw (-5.5,2pt) -- (-5.5,-2pt) node[below=.01] {$-2s-r$};
    \draw (-7,2pt) -- (-7,-2pt) node[below left=.01 and -.6] {$-2s-2r$};
    \draw (2pt,1) -- (-2pt,1) node[right] {$2t$};
    \draw (2pt,2) -- (-2pt,2) node[right] {$2t+r-1$};
    \draw (2pt,3) -- (-2pt,3) node[right] {$2t+2r-2$};
    \fill (-5.5,2) circle (2pt);
    \draw [->,very thick] (-4,1) -- (-1,1) -- (-1,5);
    \draw [->,very thick] (-7,3) -- (-2.5,3) -- (-2.5,5);
\end{tikzpicture} 
}
\caption{$E^2(  \e^\infty (D_{rst}))$}\label{F:e2orig}
\end{figure}
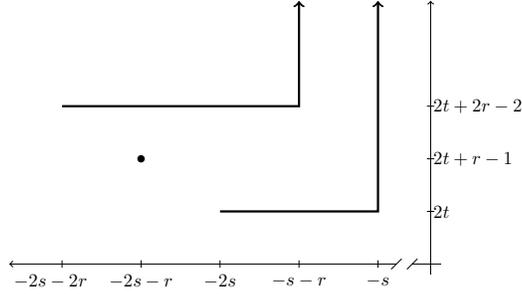

Let us first recall the construction from \cite{me1}. There, we were working in the case $n=\infty$, and the $E_2$ and $E_r$ pages of $\e^\infty(D_{rst})$ were equal and took the form of Figure~\ref{F:e2orig}. 
The idea was that the  the lower `\reflectbox{L}' in the spectral sequence for $\e^\infty(D_{rst})$ should map to the external operations on $[y]$.

We can mimic this construction here, instead using Figure~\ref{F:e2trunc}.
We first define operations on $r$-cycles. Recall that $\dmap_y$ is the universal map given in Proposition~\ref{P:universalprop}.
\begin{definition} Let
\[ \preext^m: Z^r_{-s,t}(Y) \to E^r(\e^n(Y)) \] 
be the functions defined, for $y\in Z^r_{-s,t}(Y)$, by
\begin{align*} \preext_v^{m}(y) &= E^r(\e^n(\dmap_y)) (\exte_{-s,m+t}) & m&\in  [t,t-s+n]\\
\preext_h^{m}(y) &= E^r(\e^n(\dmap_y)) (\exte_{m-s-t,2t}) &  m&\in[t-s,\min(t,t-s+n)]
 \end{align*}
where  $\exte_{p,q} \in E^2_{p,q}(\e^n(D_{rst}))$ is nonzero.
\end{definition}
Unfortunately, this definition does not make any sense yet. Why should the indicated $\exte$ survive to page $r$? How do we know that there is only one generator in the indicated bidegrees? This section is devoted to these questions.

We expect that the top operation will not be additive, so we cannot immediately carry out the program given in section~\Sopdefinf~of \cite{me1} to induce operations on the spectral sequence. We will define the Browder operation in section~\ref{S:browder} in order to study the deviation from additivity of $\preext^{t-s+n}$. 
Proposition~\ref{P:overlaps} below shows that the $\preext^m$ are well-defined, while the  computations which follow will be used in section~\ref{S:truncdef} to show that the $\preext^m$ induce operations at the spectral sequence (rather than cycle) level.
In particular, we compute enough of the differentials of the spectral sequence associated to $\e^n(D_{rst})$ to give a partial analogue to \cite[\corvanishingofbidegrees]{me1}.

The main tool is the comparison
\[ \phi: \e^n (D_{rst}) \to \e^\infty (D_{rst}) \] induced by the inclusion
\[ \sk_n W \hookrightarrow W. \]

\begin{proposition}\label{P:e2phi} Let $n\geq 1$ and $\infty > r \geq 2$. The kernel of $\phi$ on the second page is
\[ \ker (E^2 (\phi)) = \K(\rft \sqcup \drt \sqcup \rfmtwo). \]
\end{proposition}
\begin{proof} Note that the map $C(\phi)$ is just an inclusion. 
At $E^1$ the representatives of $\wrft, \wdrt$, and $\wrfmtwo$ are all vertical boundaries, so $\K (\rft \sqcup \drt \sqcup \rfmtwo) \ci \ker(E^2(\phi))$. Comparing representatives in \cite[\thmeonebasis]{me1} with the representatives in Theorem~\ref{Ts:e1basis} gives that this inclusion is equality.
\end{proof}

Define a set of integral lattice points $L=L_{rst}^n$ by
\[ [-2s-2r, -2s-2r+n]\times \set{2t+2r-2} \]
if $n\leq s+r-1$ and by
\begin{multline*} \Big( [-2s-2r,-s-r-1] \times \set{2t+2r-2} \Big) \\ \cup \Big( \set{-s-r} \times [2t+2r-2,2t+r-s-2+n] \Big) \end{multline*}
if $n\geq r+s$. 

\begin{proposition}\label{P:overlaps}
If $(p,q)\in L$, then
\[ E_{p,q}^2(\e^n (D_{rst})) = \K. \]
\end{proposition}
\begin{proof}
We know that $1$ is a lower bound for dimension since at each of these lattice points there is an element of $\drl$. The only classes at $E^2$ which might share a bidegree with $\drl$ (and hence with $L$) are $\rft$ and $\rfmtwo$. Notice that the lattice points of $L$ cover the following range of total degrees:\[ [2t-2s-2, 2t-2s+n-2], \] while by Theorem~\ref{Ts:e2page} we know that $\rfmtwo$ lives in total degree $2t-2s+n-1$ and $\rft$ lives in total degree $2t-2s+n$ and above.
\end{proof}

\begin{lemma}\label{L:emmone} Consider $E^2(\e^n(D_{rst}))$ and let 
\begin{align*}\set{c} &= \rfmone_{-2s-r,2t+r-1} \\ \set{v} &= \drb_{-2s-2r,2t+2r-2}.\end{align*} Then
\[ \ssd^r[c] = [v] \neq 0. \]
\end{lemma}
\begin{proof}
Assume that $t=0$. 
We list the ranges of total degrees of each of the various subsets which constitute a basis for $E^2$ (when $s>0$):
\begin{align*}
\rfl &: [-2s,-s+n-1] & \drl &: [-2s-2,-s+r+n-3] \\
\rft &: [-2s+n,-s+n-2] & \drt &: [-2s+n-2, -s+r+n-4] \\
\rfmone &: \set{-2s-1} & \rfmtwo &: \set{-2s+n-1}
\end{align*}
Then element $v$ is in the smallest possible total degree $-2s-2$ so must be hit by something in total degree $-2s-1$ since $E^\infty=0$ (Proposition~\ref{Ps:einfzero}). The only elements which are in total degree $-2s-1$ are $c$ and, if $n=1$, the element in $\drt_{-2s-2r,2r-1}$ (here we use that $n\geq 1$).  Examination of bidegrees indicates that the second of these could only hit $v$ via $\ssd^1$, so we have the stated result for $s>0$. The proof for $s=0$ is similar.
\end{proof}

For convenience, we recall the following computation of all nontrivial differentials in the spectral sequence associated to $\e^\infty(D_{rst})$ from \cite[\propnontrivialdiff]{me1}.
\begin{proposition}\label{P:nontrivialdiff}
The following differentials in the spectral sequence associated to $\e^\infty(D_{rst})$ are nontrivial:
\begin{align*}
\ssd^r: E^r_{-2s-r,2t+r-1} &\to E^r_{-2s-2r,2t+2r-2} \\
\ssd^{2r-1}: E^{2r-1}_{p,2t} &\to E^{2r-1}_{p-2r+1,2t+2r-2} & p&\in[-2s,-s-1] \\
\ssd^{2r-1-b}: E^{2r-1-b}_{-s,2t+b} &\to E^{2r-1-b}_{b+1-2r-s,2t+2r-2} & b&\in[0,r-2] \\
\ssd^r: E^r_{-s,2t+b} &\to E^{r}_{-s-r,2t+b+r-1} & b&\in [r-1,\infty)
\end{align*}
\end{proposition}

\begin{proposition}\label{P:diffsame}
Suppose that $v\in \rfl$ has total degree in \[ [2t-2s,2t-2s-1+n] \] and $j$ is such that $\ssd^j[\phi v]_j \neq 0$. Then 
\[ 0 \neq \ssd^j[v]_j \in \drl. \]
\end{proposition}
\begin{proof} Assume $t=0$. If $[v]_j$ makes sense (that is $\ssd^k [v]_k = 0$ for $k<j$), then
\[ \phi \ssd^j [v]_j = \ssd^j \phi [v]_j = \ssd^j [\phi v]_j \neq 0\]  so $\ssd^j [v]_j \neq 0$. Proposition~\ref{P:overlaps} coupled with Proposition~\ref{P:nontrivialdiff} then tell us that it must land in the stated place.

We now show that $\ssd^k [v]_k = 0$ for $2\leq k < j$. By Lemma~\ref{L:emmone} we know that $v$ does not hit $\rfmone$ nontrivially. 
The differential of $v$ is in the following range of total degrees,
\[ [-2s-1, -2s+n-2]\]
so we see (as in the proof of Proposition~\ref{P:overlaps}) that $v$ cannot hit any of the bidegrees spanned by $\rft$ or $\rfmtwo$. On the other hand, $\drt$ lives in the following range of total degrees
\[ [-2s+n-2, r-s+n-4], \] so it's \emph{possible} that $v$ hits something in $\drt$ if it has total degree $-2s-1+n$. But $\drt$ is so far away that this must happen on a page bigger than $j$ (see Figure~\ref{F:e2trunc} on page~\pageref{F:e2trunc}). To be precise, the differential would be one of (writing $c\in \rft_{-2s-2r,2r+n-2}$)
\begin{align*}
\ssd^{2r+s}[v] &= [c] & \text{ if } v\in \rfc_{-s,-s-1+n}  \\
\ssd^{2r+n-1}[v] &= [c] & \text{ if } v\in  \rfb_{-2s-1+n,0} 
\end{align*}
whereas $j \leq 2r-1$ by Proposition~\ref{P:nontrivialdiff}.
\end{proof}

This proposition tells us that the spectral sequence for $\e^n(D_{rst})$ vanishes in the bidegrees of $L$ at the same time as in $\e^\infty (D_{rst})$. This is precisely what we will need to help us show that the $\preext^m$ vanish on appropriate boundaries.

\section{Additivity and the Browder operation}\label{S:browder}

We would like to mimic section~\Sopdefinf~of \cite{me1}, but on first glance it appears that \cite[\propadditivity]{me1} fails when $m=t-s+n$ because of the classical formula \cite[6.5]{may}
\[ \xi_n (x+y) = \xi_n(x) + \xi_n(y) + \lambda_n (x,y) \]
where $\lambda_n$ is the Browder operation and $\xi_n (x_q) = Q^{q+n}(x)$. Surprisingly, additivity holds for the top operation as long as $s>0$. We will see in a moment that this happens because the Browder operation lands in a lower filtration degree, but first we prove the additivity statement.

\begin{proposition}[Additivity]\label{P:finiteadditivity} Let $r\geq 2$ and
\[ b= \begin{cases} t-s+n & s>0 \\ t+n-1 &s=0. \end{cases} \]
The functions
\begin{align*}
\preext_v^m:& Z_{-s,t}^r(Y) \to E^r_{-s,m+t}(\e^n(Y)) &  m&\in[t,b] \\
\preext_h^m:& Z_{-s,t}^r(Y) \to E^r_{m-s-t,2t}(\e^n(Y)) & m&\in[t-s,\min (t,t-s+n)]
\end{align*} are homomorphisms.
\end{proposition}
\begin{proof}
As in the proof of \cite[\propadditivity]{me1}, we have
\[ E^1 (\sk_n W \ten D_{rst} \ten D_{rst}) \cong H_*(\sk_n W) \ten E^1(D_{rst}) \ten E^1(D_{rst}) \] 
which is nonzero only in the following list of bidegrees:
\begin{gather*}
(-2s,2t), (-2s-r, 2t+r-1), (-2s-2r,2t+2r-2) \\
(-2s,2t+n),  (-2s-r, 2t+r-1+n) , (-2s-2r, 2t+2r-2+n)
\end{gather*}
since $H_*(\sk_n W) = \K e_0 \oplus \K (1+\sigma) e_n$. The only possible overlap with bidegrees of the operations are $(-2s,2t)$, which is the external square, and, if $s=0$, $(0,2t+n)$. But this last bidegree corresponds to the operation $\preext_v^{t+n}$ (when $s=0$), which the one operation that is excluded from the statement of the proposition.
\end{proof}

\begin{definition}[Browder Operation]\label{D:bo} Let $Y$ be a cosimplicial chain complex.
Consider $\K$ as a chain complex in degree $0$. Using the map $ \K \to \Sigma^{-n}\sk_n W$ which sends $1$ to $(1+\sigma)e_n$ for the middle arrow below, we consider the map of bicomplexes
\[ \xymatrix{   C(Y)\ten C(Y)  \ar@{->}[r]^-{AW} &  C(Y\ten Y) \ar@{->}[r] & \Sigma^{-n} C(\sk_n W \ten Y \ten Y) \ar@{->}[d] \\ && \Sigma^{-n} C(\e^n(Y))
}\]
Then we have a map
\begin{equation*} \extbrow_n: E^r_{-s,t}(Y) \ten E^r_{-s',t'}(Y) \to E^r_{-s-s',t+t'+n} (\e^n(Y)) \end{equation*} which we call the \emph{external Browder operation}.
\end{definition}

There is a discrepancy in bidegrees. Our top operation for an element in $Z^r_{-s,t}$ is in bidegree $(-s,2t+n-s)$ or $(-2s+n,2t)$, whereas the Browder operation of two elements in $E^r_{-s,t}$ is in bidegree $(-2s,2t+n)$. When $s=0$ the Browder operation still measures the deviation from additivity of the top operation.

\begin{proposition}\label{P:browderadd}
Suppose that $x,y\in Z^r_{0,t}(Y)$. Then
\[ \preext_v^{t+n}(x+y) = \preext_v^{t+n} (x) + \preext_v^{t+n}(y) + \extbrow_n([x]_r, [y]_r). \]
\end{proposition}

\noindent Before beginning the proof, we make note of the following variation of \cite[\lemsumsplit]{me1}:
\begin{lemma}\label{L:sumsplitting}
Let $X$ and $Y$ be cosimplicial chain complexes. Then
\[ \e^n( X\oplus Y)\cong \e^n(X) \oplus \e^n(Y) \oplus (\sk_n W\ten X \ten Y) \]
via
\[ e_i \ten (x+y) \ten (x'+y') \mapsto \begin{gathered} e_i \ten x \ten x' + e_i\ten y\ten y' \\ + e_i \ten x \ten y' + \sigma e_i \ten x' \ten y. \end{gathered} \]
\end{lemma}

\begin{proof}[Proof of Proposition~\ref{P:browderadd}]
Examine the diagram from \cite[\propadditivity]{me1}
\[ \xymatrix{
D_{r0t} \ar@{->}[dr]_{\dmap_{x+y}} \ar@{->}[rr] && D_{r0t} \oplus D_{r0t} \ar@{->}[dl]^{\dmap_x \oplus \dmap_y} \\
& Y 
}\] where the top map is the diagonal. 
According to Lemma~\ref{L:sumsplitting}, we have the decomposition
\[ \xymatrix{
\e^n( D) \ar@{->}[dr]_{\e^n(\dmap_{x+y})} \ar@{->}[rr] && \e^n(D \oplus D) \ar@{->}[dl] \ar@{=}[r]^-{\ref{L:sumsplitting}} & \e^n(D) \oplus \e^n(D) \oplus (\sk_n W\ten D \ten D) \ar@/^1.5pc/[dll]|{\e^n(\dmap_x) + \e^n(\dmap_y) + ???} \\
& \e^n(Y) 
}\]
where $D=D_{r0t}$. The image of $\exte_{0,2t+n}$ under $\e^n(\dmap_x)$ and $\e^n(\dmap_y)$ give $\preext^{t+n}_v(x)$ and $\preext^{t+n}_v(y)$.
We now seek to identify the image of $\exte_{0,2t+n}$ under the composite
\[ \e^n(D) \to \sk_n W \ten D \ten D \to \e^n (Y).\]

For maps $f:A\to C$ and $g:B\to C$, the following commutes
\[ \xymatrix{ \sk_n W \ten A \ten B \ar@{->}[r]^{1\ten f \ten g} \ar@{->}[d] & \sk_n W \ten C \ten C \ar@{->}[d] \\
\e^n(A\oplus B) \ar@{->}[r]^-{\e^n (f + g)} & \e^n (C)
}\]
where the left vertical arrow is the inclusion from Lemma~\ref{L:sumsplitting}. Replacing $A=B=D_{r0t}$ and $C=Y$, we extend this to the diagram
\[ \xymatrix@C=2cm@R=.5cm{ 
 C(D) \ten C(D) \ar@{->}[r]^{C(\dmap_x) \ten C(\dmap_y) } \ar@{->}[d]^{AW} &  C(Y) \ten C(Y) \ar@{->}[d]^{AW} \\
 C(D\ten D) \ar@{->}[r]^{C(\dmap_x \ten \dmap_y)} \ar@{->}[d] &  C(Y\ten Y) \ar@{->}[d] \\
\Sigma^{-n} C(\sk_n W \ten D \ten D) \ar@{->}[r]^{1\ten \dmap_x \ten \dmap_y} \ar@{->}[d] & \Sigma^{-n} C(\sk_n W \ten Y \ten Y) \ar@{->}[d] \\
\Sigma^{-n} C(\e^n(D\oplus D)) \ar@{->}[r]^-{\e^n (\dmap_x + \dmap_y)} & \Sigma^{-n} C(\e^n (Y))
}\]
The composite of the vertical maps on the right is what was used to define the external Browder operation, so
\[ \xymatrix{
E^r(D)\ten E^r(D) \ar@{->}[r] & E^r(Y)\ten E^r(Y) \ar@{->}[d] \\
& E^r(\e^n(Y)) 
}\] takes $\imath \ten \imath$ to $\extbrow_n([x]_r, [y]_r)$.
Furthermore, the Alexander-Whitney map is particularly simple on elements in cosimplicial degree $0$: $AW(\id_{[0]} \ten \id_{[0]}) = \id_{[0]} \ten \id_{[0]}$. So the vertical maps on the left give 
\[ \xymatrix@R=.2cm{
 C(D) \ten C(D) \ar@{->}[r] &  C(D\ten D) \ar@{->}[r] & \Sigma^{-n} C(\sk_n W \ten D \ten D) \\
\id_{[0]} \ten \id_{[0]} \ar@{|->}[r]  & \id_{[0]} \ten \id_{[0]} \ar@{|->}[r] & (1+\sigma)e_n \ten \id_{[0]} \ten \id_{[0]}.
}\]
At $E^1$ this coincides with the image of $\exte_{0,2t+n}$ by Lemma~\ref{L:forbrowderadd}. 
\end{proof}

\begin{lemma}\label{L:forbrowderadd}
Let $C$ be a chain complex. Consider the composite
\[ \e^n(C) \overset{\e^n \Delta}{\longrightarrow} \e^n(C\oplus C) \twoheadrightarrow \sk_n W \ten C \ten C \]
where the projection map is the one given by Lemma~\ref{L:sumsplitting}:
\begin{align*}
\e^n (X\oplus Y) &\to \sk_n W \ten X \ten Y \\
e_m\ten (x+y) \ten (x'+y') &\mapsto e_m \ten x \ten y' + \sigma e_m \ten x' \ten y.
\end{align*} Then the homology of the composite
sends $e_n \ten [c] \ten [c]$ to $(1+\sigma)e_n \ten [c] \ten [c]$.
\end{lemma}
\begin{proof} Fix a quasi-isomorphism $C\to HC$.
The following commutes, 
\[ \xymatrix{
H_* \e^n(C) \ar@{->}[d]^\cong \ar@{->}[r] & H_*(\e^n(C\oplus C)) \ar@{->}[d]^\cong \ar@{->}[r] & H_* (\sk_n W \ten C \ten C) \ar@{->}[d]^\cong \\
H_* \e^n(HC)\ar@{->}[r] &H^* \e^n (HC \oplus HC) \ar@{->}[r]&  H_*(\sk_n W \ten HC \ten HC )
}\]
so it is enough to prove that for a module $M$,
\[ H_*(\e^n (M)) \to H_*(\sk_n W  \ten M \ten M ) \] 
sends $e_n \ten m \ten m$ to $(1+\sigma) e_n \ten m \ten m$. This is an easy computation.
\end{proof}

\begin{remark} 
The formula given in Proposition~\ref{P:browderadd} says that if $y$ happens to be in $B_{0,t}^r$ then
\[ \preext^{t+n}_v (x+y) = \preext^{t+n}_v(x) + \preext^{t+n}_v(y)\] since $[y]_r=0$.  So if we show that $\preext^{t+n}_v(y)=0$ for $y\in B_{0,t}^r$ then we will know that $\preext^{t+n}_v$ induces a \emph{function} \[ E^r_{0,t} (Y) \to E^r_{0,2t+n} (\e^n (Y)). \]
\end{remark}

\section{Definition of operations}\label{S:truncdef}

The proofs of nearly everything in section~\Sopdefinf~of \cite{me1} now go through, with perhaps the only subtle point that the analogue of \cite[\lemverticalpartial]{me1} relies on the vanishing statement Proposition~\ref{P:diffsame}.

\begin{lemma}\label{L:n.lowerfiltr}
The homomorphisms $\preext^m$ vanish on $Z^{r-1}_{-s-1,t+1}$ for $r\geq 2$.
\end{lemma}
\begin{proof}
Write $r'=r-1, s'=s+1, t'=t+1$ and let $y\in Z^{r'}_{-s',t'}(Y) \ci Z^r_{s,t}(Y)$. 
Then the following commutes
\[ \xymatrix{
D_{rst} \ar@{->}[dr]_{\dmap_{y}^r} \ar@{->}[rr]^{\dmap_{\imath}} && D_{r's't'} \ar@{->}[dl]^{\dmap_y^{r'}} \\
& X 
}\]
If $r'\geq 2$, then Theorem~\ref{Ts:e2page} says that $E^r(\e^n(D_{r's't'}))$ is zero on the ranges $\set{-s}\times [2t,2t+n-s]$ and $[-2s,\min(-s-1,2t+n-s)]\times \set{2t}$ we are interested in. The diagram
\[ \xymatrix{
E^r(\e^n(D_{rst})) \ar@{->}[dr] \ar@{->}[rr] && E^r(\e^n(D_{r's't'})) \ar@{->}[dl] \\
& E^r(\e^n(Y)) 
}\]
commutes and the rightmost composition takes $\exte_{p,q}$ to zero for $(p,q)$ in the appropriate range, so all of the $\preext$ must vanish on $x$.
If $r=2$ then $E^2(\e^n(D_{1s't'})) = 0$.
\end{proof}

In particular, this shows that the $\preext^m$ vanish on $\partial F^{-s-1}$, and the proof of the following is a minor variation of that of the corresponding proposition in section~\Sopdefinf~of \cite{me1}.

\begin{proposition}
The homomorphisms $\preext^m$ vanish on $\partial F^{-s}$. \qed
\end{proposition}

\begin{lemma}\label{L:n.verticalpartial}
The vertical maps $\preext_v$ vanish on $\partial Z_{-s+r-1,t-r+2}^{r-1}$ for $r>2$.
\end{lemma}
\begin{proof} Notice that we may assume that $n\geq s$, otherwise we have not defined the vertical maps and the statement is vacuously true.
Let $r'=r-1$, $s'=s+1-r$, $t'=t-r+2$. We may assume that $y\in Z_{-s+r-1,t-r+2}^{r-1}$ has the form \[ y = \sum_{j=s-r+1}^{s-1} y^j. \]  The following diagram commutes
\[ \xymatrix{
D_{rst} \ar@{->}[dr]_{\dmap_{\partial y}} \ar@{->}[rr]^{\dmap_{\partial \imath}} && D_{r's't'} \ar@{->}[dl]^{\dmap_y} \\
& Y 
}\]
Applying Proposition~\ref{P:overlaps} to $(r',s',t')$, we see that 
the vector space $E^2_{p,q}(\e^n (D_{r's't'}))$ is one-dimensional for $p=-s$ and $q\in [2t,2t-s+n]$. Furthermore, Proposition~\ref{P:diffsame} tell us that all of these classes vanish at page $r'+1 = r$. These are exactly the bidegrees where we have defined vertical operations, so applying $E^r(\e^n(-))$ to the above diagram we see that $\preext_v(\partial y)=0$ on $E^r$.
\end{proof}

\begin{lemma}
If $r=2$ then the homomorphisms $\preext$ vanish on $\partial Z_{-s+1,t}^{1}$.
\end{lemma}
\begin{proof} As in \cite{me1}.
\end{proof}

\begin{theorem}\label{T:n.inftymain} The maps above define functions
\begin{align*}
\ext_v^m:& E_{-s,t}^r(Y) \to E^r_{-s,m+t}(\e^n(Y)) &  m&\in[t,t-s+n] \\
\ext_h^m:& E_{-s,t}^r(Y) \to E^w_{m-s-t,2t}(\e^n(Y)) & m&\in[t-s,\min(t,t-s+n)]
\end{align*}
where \[ w = \begin{cases} r & m=t-s \\ 2r-2 & m\in [t-s+1,t-r+2] \\ r+t-m & m\in [t-r+3, t]. \end{cases} \]
They are homomorphisms unless $s=0$ and $m=t-s+n$, in which case there is an error term given by Proposition~\ref{P:browderadd}.
\end{theorem}
\begin{proof}
The only missing ingredient is the vanishing of $\preext_h^m$ on an element $\partial y$ where 
$y\in Z_{-s+r-1,t-r+2}^{r-1}$ is of the form \[ y = \sum_{j=s-r+1}^{s-1} y^j. \] This is an extension of the proof of Lemma~\ref{L:n.verticalpartial}. According to Propositions~\ref{P:overlaps} and \ref{P:diffsame} applied to $(r',s',t')$, part of \cite[\corvanishingofbidegrees]{me1} applies in the spectral sequence for $\e^n(D_{r's't'})$ to give appropriate vanishing in the range of bidegrees  \[ [-2s+1,\min(-2s+n, -s-1)] \times \set{2t}.\] 
In particular,
\[ E^{2r-2}_{p,2t}(\e^n(D_{r's't'}))=0 \] 
for $p\in [-2s+1,-r-s+2] \cap I$ and
\[ E^{r-s-p}_{p,2t}(\e^n(D_{r's't'}))=0 \] 
for $p\in [-r-s+3,-s] \cap I$, where $I= [-2s+1,\min(-2s+n, -s-1)]$. 
Furthermore, Lemma~\ref{L:emmone} tells us that 
\[ E^r_{-2s,2t}(\e^n(D_{r's't'})) =0.\] 
The statement then follows by going from $p$ to $m=t+s+p$.
\end{proof}

\newcommand{\inthom}{\mathbf{hom}}
\renewcommand{\inthom}{\operatorname{map}}
\newcommand{\sethom}{\hom}
\newcommand{\wild}[1]{\widetilde{#1}}
\newcommand{\bq}{{\bar{q}}}
\newcommand{\bk}{{\bar{k}}}
\newcommand{\bj}{{\bar{j}}}
\newcommand{\shuff}[1]{\mathcal{S}_{#1}}

\newcommand{\fibbk}{\widehat{\bkspace}}
\newcommand{\zh}{Z_{h\pi}}
\newcommand{\bw}{\overline{W}}
\newcommand{\coact}{\mathbf{c}}
\newcommand{\slush}[1]{\underline{#1}}

\newcommand{\inta}{\mathbf{\rho}}
\newcommand{\intb}{\mathbf{\chi}}

\newcommand{\faked}{a}

\newcommand{\mess}{\ell}

\section{Convergence of operations}\label{S:finiteconv}

In this section, we take ``space'' to mean ``simplicial set,'' as we did in \cite{me2a}.  If $X$ is a cosimplicial space, then the target for the homology spectral sequence is $H_*(\Tot X)$, using a filtration given by Bousfield in \cite{bousfield} that we now briefly recall. 
If $V$ is a cosimplicial simplicial module and $\mess\leq \infty$, then there is a natural quasi-isomorphism 
\[ \xymatrix@R=7pt{ 
N \Tot_\mess V \ar@{->}[r]^{\phi_\mess}  \ar@{=}[d] & 
T_\mess CN V \ar@{=}[d] \\
N \inthom (\sk_\mess \Delta, V) &
 TCNV / F^{-(\mess + 1)}
}\]  
by \cite[2.2]{bousfield}.
The composite
$S_* \Tot X = N \K \Tot X \to N \Tot (\K X) \overset{\phi_\infty}{\to} T CN \K X$ 
is used to define the filtration on $H_*(\Tot X)$, which is given by
\[ F^{-s} = \ker(H \Tot X \to H TCN\K X \to H T_{s-1} CN\K X ),  \]
so $F^{-(s+1)} \ci F^{-s}$. This gives the abutment map
\begin{equation*}  (F^{-s} / F^{-s-1}) H (\Tot X) \hookrightarrow E_{-s}^\infty (X; \K). \label{SYMBOL:ABUTMENT} \end{equation*}

The arguments of \cite{me2a} may be modified to give
\begin{theorem}\label{T:convergencefinite} Suppose that $X$ is a cosimplicial space and $m\leq t-s+n$. Then
\[ Q^m [ F^{-s} H_{t-s} (\Tot (X))] \ci F^{-v}H_{t-s+m} (\Tot(S^n \times_\pi X^{\times 2})) \]
where \[ v(t,s,m) = \begin{cases} s & m\geq t \\ t+s-m & t-s \leq m \leq t. \end{cases} \]Furthermore, for each $m$  the following diagram commutes:
\[ \xymatrix{
(F^{-s} / F^{-s-1}) H_{t-s} (\Tot X) \ar@{->}[r] \ar@{->}[dd]^{Q^m} & E_{-s,t}^\infty (S_*(X)) \ar@{->}[d]^{\ext^m} \\
& E_{-v, v-s+m+t}^\infty (\sk_n W\tp S_*(X)^{\ten 2}) \ar@{->}[d]^{\cong}_{E^1 - iso} \\
(F^{-v} / F^{-v-1}) H_{t-s+m} (\Tot (S^n \times_\pi X^{\times 2})) \ar@{->}[r] & E_{-v, v-s+m+t}^\infty (S_*(S^n \times_\pi X^{\times 2}))
}\]
\end{theorem}
The obvious necessary change is to replace, throughout \cite{me2a}, $E\pi$ and $B\pi$ by $S^n = \sk_n E\pi$ and $\mathbb{R} P^n = \sk_n B\pi$, respectively.
The one other change we must make in comparison with the proof of the analogous result from \cite{me2a} is that we must use the results of section~\ref{S:truncinf} to get a handle on the $E^\infty$ term for the spectral sequence of $\e^n(D_{\infty st})$.

In particular, Theorem~\ref{T:convergencefinite} implies convergence of the internal operations when $X$ is a cosimplicial $\msc_{n+1}$-space:

\begin{theorem}
Suppose that $X$ is a cosimplicial $\msc_{n+1}$-space. Then the Araki-Kudo operation $Q^m$ has the following effect on the filtration:
\[ Q^m [ F^{-s} H_{t-s} (\Tot (X))] \ci F^{-v}H_{t-s+m} (\Tot X), \]
where $v$ is given above. 
Furthermore, for each $m \in [t-s,t-s+n]$,  the following diagram commutes:
\[ \xymatrix{
\left({F^{-s}}/{F^{-s-1}}\right) H_{t-s} (\Tot X) \ar@{->}[r] \ar@{->}[d]^{Q^m} & E_{-s,t}^\infty (X) \ar@{->}[d]^{\ext^m}\\
\left( {F^{-v}}/{F^{-v-1}}\right) H_{t-s+m} (\Tot (X)) \ar@{->}[r] & E_{-v, v-s+m+t}^\infty (X)
}\]
\end{theorem}

We also note that the results of \cite{me2a} which are related to multiplication, namely 
Theorem 4.4 and Corollary 4.5, may be extended to the present situation with only minor changes to the proofs that appear in \cite{me2a}.

\begin{theorem}\label{T:fmult} Fix $n\geq 2$, and let $X$ be a cosimplicial space. Then the external multiplication $\mu: H \Tot X \ten H \Tot X \to H \Tot (S^n \times_\pi X^{\times 2})$ is compatible with the filtration in the sense that 
\[\mu(F^{-s} \ten F^{-s'}) \ci F^{-s-s'}.\] Furthermore, we have that the following diagram commutes.
\[ \xymatrix@C-10pt{
(F^{-s} / F^{-s-1}) H\Tot X \ten (F^{-s'} / F^{-s'-1}) H\Tot X \ar@{->}[r] \ar@{->}[d]   & E^\infty_{-s} (S_*(X)) \ten  E^\infty_{-s'} (S_*(X)) \ar@{->}[d] \\
(F^{-s-s'}/F^{-s-s'-1}) H\Tot (S^n \times_\pi X^{\times 2}) \ar@{->}[r] & E^\infty_{-s-s'} (S_*(S^n \times_\pi X^{\times 2})) \\
}\] If  $X$ is a cosimplicial $\msc_{n+1}$-space, then we have the corresponding internal statement that the diagram
\[ \xymatrix@C-10pt{
(F^{-s} / F^{-s-1}) H\Tot X \ten (F^{-s'} / F^{-s'-1}) H\Tot X \ar@{->}[r] \ar@{->}[d]   & E^\infty_{-s} (S_*(X)) \ten  E^\infty_{-s'} (S_*(X)) \ar@{->}[d] \\
(F^{-s-s'}/F^{-s-s'-1}) H\Tot X \ar@{->}[r] & E^\infty_{-s-s'} (S_*(X)) \\
}\]
commutes.  \qed
\end{theorem}

We now turn to the Browder operation, which is fundamentally an algebraic operation. Recall that in \cite{me2a} we proved a theorem which stated that in homology the map $\phi_\mess$ is compatible with tensor products, in the sense that the diagram
\[ \xymatrix{
H(\Tot_\mess U) \ten H(\Tot_\mess V) \ar@{->}[r]^{\phi_\mess \otimes \phi_\mess} \ar@{->}[d]  &  H (T_\mess CNU) \ten H (T_\mess CNV) \ar@{->}[d] \\
H (\Tot_\mess (U\ten V)) \ar@{->}[r]^{\phi_\mess} & H(T_\mess CN(U\ten V))
}\] commutes.
Regarding $\K S^n$ also as a constant cosimplicial simplicial $\K$-module, we can iterate this theorem to get that bottom of the diagram
\[ \xymatrix{ 
H_*(\Tot_\mess U) \ten H_*(\Tot_\mess V) \ar@{}[dr]|{\#} \ar@{->}[d]^{(1+\sigma) e_n \ten - \ten -} \ar@{->}[r]  &H_* (T_\mess CNU) \ten H_* (T_\mess CNV) \ar@{->}[d]^{(1+\sigma) e_n \ten - \ten -} \\
H_*(S^n) \ten H_*(\Tot_\mess U) \ten H_*(\Tot_\mess V) \ar@{->}[r] \ar@{->}[d]  & H_*(S^n) \ten H_* (T_\mess CNU) \ten H_* (T_\mess CNV) \ar@{->}[d] \\
H_* (\Tot_\mess (\K S^n \ten U\ten V)) \ar@{->}[r] & H_*(T_\mess CN(\K S^n \ten U\ten V))
}\]
commutes. The left hand composite 
\begin{multline*} N\Tot_\mess U \ten N\Tot_\mess V \to \sk_n W \ten N\Tot_\mess U \ten N\Tot_\mess V \\ \overset{\nabla}{\to} N[ \K S^n \ten \Tot_\mess U \ten \Tot_\mess V ] \to N \Tot_\mess (\K S^n \ten U \ten V) \end{multline*}
can be seen a cosimplicial extension of the usual external Browder operation (see \cite[p.184]{may}).
The right hand composite is given by
\begin{align*} T_\mess CN U \ten T_\mess CN V & \to \sk_n W \ten T_\mess CN U \ten T_\mess CN V \\
&\to T_\mess ( (\sk_n W)^v \ten CN U \ten CN V ) \\
&\overset{AW}{\to} T_\mess C(\sk_n W \ten NU \ten NV) \\
&\overset{\nabla}{\to} T_\mess CN (\K S^n \ten U \ten V). 
\end{align*}
The composite \begin{multline*} T_\mess CN U \ten T_\mess CN V  \to \sk_n W \ten T_\mess CN U \ten T_\mess CN V \\ \to T_\mess ( (\sk_n W)^v \ten CN U \ten CN V ) \overset{AW}{\to} T_\mess C(\sk_n W \ten NU \ten NV) \end{multline*}
is equal to 
\begin{multline*} T_\mess CN U \ten T_\mess CN V \to T_\mess ( CN U \ten CN V ) \\ \to T_\mess [ (\sk_n W)^v \ten CN U \ten CN V ] \to T_\mess C (\sk_n W \ten NU \ten NV). \end{multline*} We see this show up as the the bottom left portion of the commutative diagram
\[ \xymatrix{
C(NU) \ten C(NV) \ar@{->}[r]^-{AW} \ar@{->}[d] & C(NU \ten NV) \ar@{->}[d] \\
(\sk_n W)^v \ten C(NU) \ten C(NV) \ar@{->}[r]^-{AW} & C(\sk_n W \ten NU \ten NV)
}\]
Our definition of the Browder operation Definition~\ref{D:bo} is essentially the upper right portion of this diagram. Thus we see that our definition agrees with the classical one.
\begin{theorem}\label{T:browderconv} We consider the map
\[ H_k(\Tot X) \ten H_{k'}(\Tot X) \to H_{k+k'+n} (S^n \times \Tot X \times \Tot X) \to H_{k+k'+n} (\Tot (S^n \times_\pi  X^{\times 2})) \] as the classical external Browder operation.
This map descends to filtration quotients and, furthermore, our spectral sequence version $\extbrow_n$ from Definition~\ref{D:bo} makes the following diagram commute.
\[ \xymatrix{
(F^{-s} / F^{-s-1}) H\Tot X \ten (F^{-s'} / F^{-s'-1}) H\Tot X \ar@{->}[r] \ar@{->}[d]   & E^\infty_{-s} (S_*(X)) \ten  E^\infty_{-s'} (S_*(X)) \ar@{->}[d]^{\extbrow_n} \\
(F^{-s-s'}/F^{-s-s'-1}) H\Tot (S^n \times_\pi X^{\times 2}) \ar@{->}[r] & E^\infty_{-s-s'} (S_*(S^n \times_\pi X^{\times 2})) \\
}\]
\end{theorem}
\begin{proof}
Most of the argument precedes the theorem statement, and the rest is as in \cite{me2a}.
\end{proof}

\ack
Parts of this work are adapted from the author's Ph.D. thesis, and the author thanks his advisor, Jim McClure, for careful readings and numerous clarifying suggestions. The exposition of this paper benefited from helpful comments from the referee.

\end{document}